\documentclass[twocolumn]{autart}    

\usepackage{graphicx}

\usepackage[utf8]{inputenc}
\usepackage{xcolor}
\usepackage{xspace}

\usepackage{amsmath,amssymb}
\usepackage{mathtools}
\usepackage{bm,fixmath} 

\usepackage{siunitx}
\usepackage{accents}
\usepackage{graphicx,tikz}
\usepackage{color}
\usepackage{bm,fixmath}
\usetikzlibrary{positioning}
\usepackage{empheq}
\usepackage{nicematrix}
\usepackage{booktabs}
\usetikzlibrary{calc,shapes,arrows.meta,decorations.pathreplacing,decorations.markings,positioning}

\usetikzlibrary{positioning}

\makeatletter

\pgfdeclareshape{record}{
\inheritsavedanchors[from={rectangle}]
\inheritbackgroundpath[from={rectangle}]
\inheritanchorborder[from={rectangle}]
\foreach \x in {center,north east,north west,north,south,south east,south west}{
\inheritanchor[from={rectangle}]{\x}
}
\foregroundpath{
\pgfpointdiff{\northeast}{\southwest}
\pgf@xa=\pgf@x \pgf@ya=\pgf@y
\northeast
\pgfpathmoveto{\pgfpoint{0}{0.33\pgf@ya}}
\pgfpathlineto{\pgfpoint{0}{-0.33\pgf@ya}}
\pgfpathmoveto{\pgfpoint{0.33\pgf@xa}{0}}
\pgfpathlineto{\pgfpoint{-0.33\pgf@xa}{0}}
\pgfpathmoveto{\pgfpointadd{\southwest}{\pgfpoint{-0.33\pgf@xa}{-0.6\pgf@ya}}}
\pgfpathlineto{\pgfpointadd{\southwest}{\pgfpoint{-0.5\pgf@xa}{-0.6\pgf@ya}}}
\pgfpathlineto{\pgfpointadd{\northeast}{\pgfpoint{-0.5\pgf@xa}{-0.6\pgf@ya}}}
\pgfpathlineto{\pgfpointadd{\northeast}{\pgfpoint{-0.33\pgf@xa}{-0.6\pgf@ya}}}
}
}
\makeatother

\tikzset{%
  saturation block/.style={%
    fill=white,
    path picture={
      \pgfpointdiff{\pgfpointanchor{path picture bounding box}{north east}}%
        {\pgfpointanchor{path picture bounding box}{south west}}
      \pgfgetlastxy\x\y
      \tikzset{x=\x*.4, y=\y*.4}
      \draw[thin,draw = gray] (-1,0) edge[{<[scale=0.5]}-] (1,0) (0,-1) edge[{<[scale=0.5]}-] (0,0.5); 
      \draw[thick] (.7,0) -- (0,0) -- (-.7,-.7) ;
    }
  }
}

\tikzset{
box/.style={
  rectangle,
  inner sep=0pt,
  minimum size=0.75cm,
  text width=1.6cm,
  align=center,
  draw=black,
  fill=white
  },
dot/.style={
   circle,inner sep=0pt,minimum size=0pt,fill=white
   }
}

\DeclareMathOperator*{\argmin}{arg\,min}

\DeclareMathOperator{\proj}{proj}

\newcommand*\widefbox[1]{\fbox{\hspace{2em}#1\hspace{2em}}}

\newcommand{\norm}[1]{\left\lVert#1\right\rVert}

\newcommand{\ubar}[1]{\underaccent{\bar}{#1}}
\renewcommand{\qed}{\hfill$\square$}

\usepackage{hhline}
\def\N{{\mathbb{N}}}
\def\R{{\mathbb{R}}}

\def\bv{{\mathbold{b}}}
\def\cv{{\mathbold{c}}}

\def\e{{\mathbold{e}}}
\def\f{{\mathbold{f}}}
\def\h{{\mathbold{h}}}

\def\vv{{\mathbold{v}}}
\def\w{{\mathbold{w}}}
\def\v{{\mathbold{v}}}
\def\x{{\mathbold{x}}}
\def\y{{\mathbold{y}}}
\def\z{{\mathbold{z}}}

\def\A{{\mathbold{A}}}
\def\B{{\mathbold{B}}}
\def\Cm{{\mathbold{C}}}

\def\F{{\mathbold{F}}}
\def\G{{\mathbold{G}}}

\def\I{{\mathbold{I}}}

\def\P{{\mathbold{P}}}
\def\Q{{\mathbold{Q}}}

\def\U{{\mathbold{U}}}
\def\V{{\mathbold{V}}}

\def\K{{\mathbold{K}}}
\def\Lm{{\mathbold{\Lambda}}}

\def\0{{\mathbf{0}}}
\def\1{{\mathbf{1}}}

\def\lmin{{\ubar{\lambda}}}
\def\lmax{{\bar{\lambda}}}

\def\Ts{{T_\mathrm{s}}}

\def\fh{{\mathbold{\hat{f}}}}
\def\hh{{\mathbold{\hat{h}}}}
\def\bh{{\mathbold{\hat{b}}}}
\def\xh{{\mathbold{\hat{x}}}}
\def\wh{{\mathbold{\hat{w}}}}
\def\eh{{\mathbold{\hat{e}}}}
\NiceMatrixOptions{cell-space-limits = 4pt}

\begin{document}

\begin{frontmatter}

\title{A Control Theoretical Approach to\\Online Constrained Optimization\thanksref{footnoteinfo}} 

\thanks[footnoteinfo]{Corresponding author N.~Bastianello.\\
This work was supported in part by the Italian Ministry of University and Research (2020–2023) through Research Project PRIN 2017 Advanced Network Control of Future Smart Grids.\\
The work of N.B. was partially supported by the European Union’s Horizon Research and Innovation Actions programme under grant agreement No. 101070162.}

\author[unipd]{Umberto Casti}\ead{umberto.casti@phd.unipd.it},
\author[kth]{Nicola Bastianello}\ead{nicolba@kth.se},
\author[unipd]{Ruggero Carli}\ead{carlirug@dei.unipd.it},
\author[unipd]{Sandro Zampieri}\ead{zampi@dei.unipd.it}

\address[unipd]{Department of Information Engineering (DEI), University of Padova, Italy}
\address[kth]{School of Electrical Engineering and Computer Science and Digital Futures, KTH Royal Institute of Technology, Sweden}

\begin{keyword}                           
online optimization, constrained optimization, control theory, anti-windup               
\end{keyword}                             

\begin{abstract}                          
In this paper we focus on the solution of online problems with time-varying, linear equality and inequality constraints. Our approach is to design a novel online algorithm by leveraging the tools of control theory. In particular, for the case of equality constraints only, using robust control we design an online algorithm with asymptotic convergence to the optimal trajectory, differently from the alternatives that achieve non-zero tracking error. When also inequality constraints are present, we show how to modify the proposed algorithm to account for the wind-up induced by the nonnegativity constraints on the dual variables. We report numerical results that corroborate the theoretical analysis, and show how the proposed approach outperforms state-of-the-art algorithms both with equality and inequality constraints.
\end{abstract}

\end{frontmatter}

\section{Introduction}\label{sec:introduction}

Recent advances in technology have made \textit{online optimization} problems increasingly relevant in a wide range of applications, \textit{e.g.} control \cite{liao-mcpherson_semismooth_2018,paternain_realtime_2019}, signal processing \cite{hall_online_2015,fosson_centralized_2021,natali_online_2021}, and machine learning \cite{shalev_online_2011,dixit_online_2019,chang_distributed_2020}.
Online problems are characterized by costs and constraints that change over time, mirroring changes in the dynamic environments that arise in these applications. The focus then shifts from solving a single problem (as in static optimization) to solving a sequence of problems \textit{in real time}.

Specifically, in this paper we are interested in online optimization problems with linear equality and inequality constraints
\begin{equation}\label{eq:quadratic-online-optimization-general}
\begin{split}
	&\x_k^* = \argmin_{\x \in \R^n} f_k(\x) \coloneqq \frac{1}{2} \x^\top \A \x + \bv_k^{\top} \x \\
	&\text{s.t. } \G \x = \h_k \quad \G^{\prime} \x \leq \h_k^{\prime}
\end{split} \qquad k \in \N
\end{equation}
where consecutive cost functions $f_k(\x)$ and sets of constraints arrive every $\Ts > 0$ seconds. In the following, we assume that each problem in~\eqref{eq:quadratic-online-optimization-general} has a unique solution.

Drawing from \cite{simonetto_timevarying_2020}, we can distinguish two approaches to online optimization: \textit{unstructured} and \textit{structured}.
Unstructured algorithms are designed by applying static optimization methods (\textit{e.g.} gradient descent-ascent) to each (static) problem in the sequence~\eqref{eq:quadratic-online-optimization-general}. The convergence of these algorithms has been analyzed in \textit{e.g.} \cite{dixit_online_2019,shalev_online_2011,hall_online_2015,fosson_centralized_2021} for unconstrained problems, and \cite{selvaratnam_numerical_2018,cao_online_2019,zhang_online_2021,lupien_online_2023} for constrained ones; see also the surveys \cite{dallanese_optimization_2020,simonetto_timevarying_2020}.
However, straightforwardly repurposing static methods for an online scenario does not leverage in any way knowledge of the dynamic nature of these problems. And indeed, usual convergence analyses guarantee that the output of the online algorithm tracks the optimal trajectory $\{ \x_k^* \}_{k \in \N}$ with only finite precision.

As a consequence, the alternative approach of structured methods has received increasing attention. ``Structured'' refers to the fact that the online algorithms are designed by incorporating information -- a model -- of the cost and constraints variability, with the goal of tracking $\{ \x_k^* \}_{k \in \N}$ with improved precision.
Structured methods applied to constrainted online problems include prediction-correction \cite{simonetto_dual_2019,bastianello_extrapolation_2023}, interior-point methods \cite{fazlyab_prediction_2018}, and virtual queue algorithms \cite{cao_virtual-queue-based_2018}. These methods leverage a set of assumptions on the rate of change of~\eqref{eq:quadratic-online-optimization-general} as a rudimentary model to provably increase tracking precision. Although they outperform unstructured methods, these algorithms cannot reach zero tracking error and constraint violation.
Recently, some of the authors have proposed a structured algorithm for unconstrained problems, which is designed using a fully fledged model of the cost function \cite{bastianello_internal_2022}. Knowledge of this model yields an improvement in the algorithm's performance.
In this paper we focus on extending this control theoretical design approach to constrained problems. Similarly to the unconstrained case, we will show that restricting to the class of problems~\eqref{eq:quadratic-online-optimization-general} it is possible to achieve zero tracking error and constraint violation.

By taking this approach, we place our contribution at the productive intersection of control theory and optimization, both static and online.
In particular, the design of optimization algorithms using control theoretical techniques has been explored in e.g. \cite{lessard_analysis_2016,scherer_optimization_2023} for static optimization, and \cite{shahrampour_online_2017,shahrampour_distributed_2018,simonetto_optimization_2023,davydov_contracting_2023} for online optimization. In particular, \cite{shahrampour_online_2017,shahrampour_distributed_2018} exploit a linear model for the evolution of $\{ \x_k^* \}_{k \in \N}$ to analyze convergence of online algorithms, \cite{simonetto_optimization_2023} leverages Kalman-filtering for \textit{stochastic} online optimization, and \cite{davydov_contracting_2023} uses contraction analysis to design and analyze continuous-time online algorithms.
In turn, online optimization has also been used extensively as a tool to design controllers, giving rise to \textit{feedback optimization}. As the name suggests, in feedback optimization an online optimization algorithm is connected in a feedback loop with a dynamical system \cite{bernstein_online_2019,colombino_online_2020,hauswirth_timescale_2021}, as in the model predictive control (MPC) set-up \cite{liao-mcpherson_semismooth_2018}. Specifically, the output of the algorithm serves as control input to the system, whose output in turn acts on the online problem formulation.

We are now ready to discuss the proposed contribution.
As mentioned above, we focus on the solution of online problems with linear equality and inequality constraints, see~\eqref{eq:quadratic-online-optimization-general}. In particular, the goal is to extend the (structured) control theoretical approach of \cite{bastianello_internal_2022} to handle these problems. The key to this approach is interpreting the online problem as the plant, for which we have a model, that needs to be controlled.
We start our development by working on quadratic online problems with linear equality constraints, which are relevant in and of themselves in signal processing and machine learning \cite{dixit_online_2019,shalev_online_2011}. This class of problems can be reformulated as \textit{linear control problems with uncertainties}, which allows us to design a novel online algorithm by leveraging the internal model principle and techniques from robust control, as well as linear algebra results on saddle matrices. The proposed algorithms can then be shown to achieve \textit{zero} tracking error, differently from unstructured methods and alternative structured ones. 
%
%
However, the proposed algorithm cannot be applied in a straightforward manner when also inequality constraints are present. Indeed, these constraints translate into non-negativity constraints on the dual variables, and the problem cannot be cast as a linear control problem. But by noticing that non-negativity constraints act as a \textit{saturation}, we are able to incorporate them in the algorithm with the use of an \textit{anti wind-up} technique. 
We conclude the paper with numerical results that compare the proposed algorithm to alternative methods, showing its promising performance when both equality and inequality constraints are present.
To summarize, we offer the following contributions: 1) We take a control theoretical approach to design online algorithms for problems with linear equality and inequality constraints. In particular, when the cost is quadratic and only equality constraints are present, we show how a zero tracking error can be achieved. 2) We extend the proposed algorithm to deal with inequality constraints as well, by leveraging an anti wind-up mechanism. 3) We present numerical simulations to show how the proposed algorithms can outperform alternative unstructured methods, both with equality and inequality constraints, and both for quadratic and non-quadratic costs.
\section{Problem formulation and background}\label{sec:background}

\noindent \textbf{Notation.} We denote by $\mathbb{N}$, $\mathbb{R}$ the natural and real numbers, respectively, and by $\mathbb{R}\left[z\right],\mathbb{R}\left(z\right)$ the set of polynomials and of rational functions in $z$ with real coefficients. Vectors and matrices are denoted by bold letters, \emph{e.g.} $\x\in\mathbb{R}^n$ and $\A \in \mathbb{R}^{n\times n}$. The symbol $\I$ denotes the identity, $\1$ denotes the vector of all $1\mathrm{s}$ and $\0$ denotes the vector of all $0\mathrm{s}$. The $\mathrm{2}$-norm of a vector and the induced $\mathrm{2}$-norm of a matrix are both denoted by $\norm{\cdot}$. Moreover, in the following, the entry-wise partial order in $\mathbb{R}^n$ is considered. Therefore, the notation $\x \leq \y$ denotes the fact that the inequality holds for every pair of corresponding entries in $\x$ and $\y$. 
The symbol $\A \preceq \B$ (or $\A \prec \B$) means that the matrix $\B-\A$ is positive semi-definite (or positive definite). $\otimes$ denotes the Kronecker product. 
For a function $f(\x)$ from $\mathbb{R}^n$ to $\mathbb{R}$, we denote by $\nabla f$ its gradient while for a function $f\left(\x,\w\right)$ from $\mathbb{R}^{n}\times\mathbb{R}^m$ to $\mathbb{R}$, we denote $\nabla_{\x}f\left(\x,\w\right)$ and $\nabla_{\w}f\left(\x,\w\right)$ its gradient with respect to $\x$ and $\w$, respectively. 
$f(\x)$ is $\lmin$-strongly convex, $\lmin > 0$, iff $f(\x) - \frac{\lmin}{2}\norm{\x}^2$ is convex, and $\lmax$-smooth iff $\nabla f(\x)$ is $\lmax$-Lipschitz continuous. 
Given a symmetric matrix $\A$, with $\lmin\left(\A\right)$ and $\lmax\left(\A\right)$ we denote the minimum and the maximum eigenvalue of $\A$, respectively, while for any matrix $\B$ with $\ubar{\sigma}\left(\B\right)$ and $\bar{\sigma}\left(\B\right)$ we denote the minimum and the maximum singular value of $\B$. Moreover with $\kappa\left(\B\right)$ we denote the condition number on $\B$. 
Finally, $\proj_{\geq 0}(\cdot): \mathbb{R}^n \to \mathbb{R}^n_{\geq 0}$ represents the function that returns the closest element (in Euclidean norm) in the non-negative orthant. Finally, we denote by $\mathcal{Z}\left[\cdot\right]$ the $\mathcal{Z}$-transform of a given signal.

\subsection{Problem formulation}\label{subsec:problem-formulation}
We consider Problem~\eqref{eq:quadratic-online-optimization-general} where $\A \in \R^{n\times n}$, $\G \in \R^{p\times n}$ and $\G^{\prime} \in \R^{p^{\prime}\times n}$ are fixed matrices, while $\bv_k \in \mathbb{R}^n$, $\h_k \in \mathbb{R}^p$ and $\h_k^{\prime} \in \mathbb{R}^{p^\prime}$ are time-varying. In addition to that, we make the following assumptions.

\begin{assum}[Strongly convex and smooth]\label{as:generic-online-problem-A}
The symmetric matrix $\A$ is such that 
\begin{equation}\label{eq:estA}
    \ubar{\nu}\I \preceq \A \preceq \bar{\nu}\I,
\end{equation} with $0<\ubar{\nu}\leq\bar{\nu} \leq + \infty$. This is equivalent to imposing that the cost functions $\{ f_k \}_{k \in \N}$ are $\ubar{\nu}$-strongly convex and $\bar{\nu}$-smooth for any time $k \in \N$.
\end{assum}

\begin{assum}[Constraints]\label{as:constraints}(i) The constraint matrix $\left[\G^{\top} \,|\,{\G^{\prime}}^{\top}\right]^{\top} \in \R^{\left(p+p^{\prime}\right) \times n}$ is full row rank. (ii) The symmetric matrix $\G\A^{-1}\G^{\top}$ is such that
\begin{equation}\label{eq:estGAGT}
    \ubar{\mu}\I \preceq \G\A^{-1}\G^{\top} \preceq \bar{\mu}\I,
\end{equation}
with $0<\ubar{\mu} \leq \bar{\mu} \leq + \infty$.
\end{assum}
Assumptions~\ref{as:generic-online-problem-A} and \ref{as:constraints}, imply that each problem in \eqref{eq:quadratic-online-optimization-general}  has a unique minimizer, and we can define the optimal trajectory $\{ \x_k^* \}_{k \in \N}$. In particular, Assumption~\ref{as:generic-online-problem-A} is widely used in online optimization \cite{dallanese_optimization_2020}, since the strong convexity of the cost implies that, at any $k\in \N$, there is a unique primal solution $\x^*$. Moreover, Assumption~\ref{as:constraints} (i) guarantees that also the dual solution is unique \cite[p. 523]{boyd_vandenberghe_2004}.
\begin{rem}\label{rem:robust}
In the following, we assume that the online algorithm only has access to an oracle for the gradient. Moreover, we assume that in the design of this algorithm we can use only the values $\ubar{\nu}$, $\bar{\nu}$, $\ubar{\mu}$ and $\bar{\mu}$ appearing in \eqref{eq:estA} and \eqref{eq:estGAGT}.
This assumption is made to align with gradient methods that our technique aims to improve, where the algorithm does not need to possess full knowledge of $\A$ and $\G$. Concerning the bound \eqref{eq:estGAGT} on $\G\A^{-1}\G^{\top}$ notice that the values of $\ubar{\mu}$ and $\bar{\mu}$ can be obtained from $\lmin$, $\lmax$ and from bounds on the singular values of $\G$. Indeed, observe that
\begin{equation}\label{eq:BoundsGAinvGT}
    \frac{\ubar{\sigma}^2\left(\G\right)}{\lmax\left(\A\right)} \leq \lmin\left(\G\A^{-1}\G^{\top}\right) \leq \lmax\left(\G\A^{-1}\G^{\top}\right) \leq \frac{\bar{\sigma}^2\left(\G\right)}{\lmin\left(\A\right)}.
\end{equation}
\end{rem}

Our objective is to design an online algorithm that can \emph{track the optimizer sequence} $\{ \x_k^* \}_{k \in \N}$ \emph{in a real-time fashion}. Let $\{ \x_k \}_{k \in \N}$ be the output of such an algorithm. The goal is to ensure that
\begin{equation}\label{eq:tracking-error}
	\limsup_{k \to \infty} \norm{\x_k - \x_k^*} \leq B < \infty.
\end{equation}
The real-time requirement on the algorithm derives from the dynamic nature of problem~\eqref{eq:quadratic-online-optimization-general}, with new cost and constraints being received every $\Ts$ seconds. Therefore, as in the broader literature \cite{simonetto_timevarying_2020,dallanese_optimization_2020}, the goal is for the online algorithm to output $\x_k$ within the time interval $[(k-1) \Ts, k \Ts)$. The standing assumption is that $\Ts$ seconds (which is often a design parameter) are enough for a gradient evaluation and computing $\x_k$ (see the discussion in \cite[sec.~II.C]{simonetto_timevarying_2020}).


\smallskip

As discussed in section~\ref{sec:introduction}, in this paper we focus on a \emph{control theoretical} approach to designing online algorithms that can achieve~\eqref{eq:tracking-error}.
We start our development by focusing on problems with equality constraints only, for which we propose a novel algorithm in section~\ref{subsec:eqConstr}, and analyze its convergence. In section~\ref{sec:ineq-constraints} we then extend the applicability of the algorithm to inequality constraints as well, by leveraging anti wind-up techniques.


\section{Online Optimization with Equality Constraints}\label{subsec:eqConstr}
In this section we will analyze Problem~\eqref{eq:quadratic-online-optimization-general} with only equality constraints, namely we will study the following problem:
\begin{equation}\label{eq:quadratic-online-optimization}
	\x_k^* = \argmin_{\x \in \R^n} f_k(\x)
	    \qquad \text{s.t.} \ \G \x = \h_k \qquad k \in \N,
\end{equation}
where recall that $f_k(\x) \coloneqq \frac{1}{2} \x^\top \A \x + \bv_k^{\top} \x$.
We also suppose that we have a partial knowledge on the time-varying terms, as clarified by the following assumption.
\begin{assum}[Models of $\bv_k$ and $\h_k$]\label{as:modelsbh}
The sequence of vectors $\left\lbrace\bv_k \right\rbrace_{k \in \N}$ and  $\left\lbrace \h_k \right\rbrace_{k \in \N}$ have rational $\mathcal{Z}$-transforms
\begin{equation}\label{eq:linearModels}
\begin{aligned}
    \bh\left(z\right) \coloneqq \mathcal{Z}\left[ \bv_k \right]= \frac{\bv_\mathrm{N} \left(z\right)}{p\left(z\right)},\ 
    \hh\left(z\right) \coloneqq \mathcal{Z}\left[ \h_k \right] = \frac{\h_\mathrm{N}\left(z\right)}{p\left(z\right)}
\end{aligned}
\end{equation}
where the polynomial $p\left(z\right) = z^m + \sum_{i = 0}^{m-1} p_i z^i$ is known and whose zeros are assumed to be all marginally stable\footnote{A zero $\bar z\in{\mathbb C}$ of a polynomial is said to be marginally stable if $|\bar z|=1$.}. The numerators, $\bv_\mathrm{N}\left(z\right)$ and $\h_\mathrm{N}\left(z\right)$, are instead assumed to be unknown.
\end{assum}


\begin{rem}\label{rem:assum}
    Observe that imposing the same polynomial at the denominator on the rational $\mathcal{Z}$-transforms in Equation~\eqref{eq:linearModels} is not restrictive. In fact, if $\bh(z)$ and $\hh(z)$ have different denominators, they can be rewritten with a common denominator. Specifically, it is always possible to multiply and divide $\bh(z)$ by a polynomial and similarly multiply and divide $\hh(z)$ by another polynomial, such that the resulting rational functions share the same denominator. This common denominator is the least common multiple of the original denominators. Consequently, the modified rational functions represent the same signal in the time domain and have the same denominator $p\left(z\right)$.
Finally, we remark that the internal model is assumed to be static, and hence $p(z)$ is a polynomial with a fixed degree independent of $k$.
\end{rem}
Taking inspiration from \cite{bastianello_internal_2022}, we develop an algorithm based on the \emph{primal-dual} problem that leverages the Internal Model Principle \cite{byrnes_output_1997} to address Problem~\eqref{eq:quadratic-online-optimization}. In this case, we are able to provide precise and rigorous proof of the algorithm convergence.

\subsection{Design of the algorithm}\label{subsec:algorithm}
In this section we extend the model-based approach proposed in  \cite{bastianello_internal_2022} to address the constrained problem~\eqref{eq:quadratic-online-optimization}. To this aim, we define the time-varying Lagrangian
\begin{equation}\label{eq:lagrangian}
	\mathcal{L}_k(\x, \w) = f_k(\x) +\w^{\top}\left(\G \x - \h_k \right), \quad k \in \N,
\end{equation}
where $\w \in\R^p$ represents the Lagrange multiplier. 
The solution of the \emph{primal-dual} problem \cite[p.~244]{boyd_vandenberghe_2004} is given by the pair of vectors $\x_k^* \in \R^n$ and $\w_k^* \in \R^p$ satisfying the equation
\begin{equation}\label{eq:error_zero}
	\begin{bmatrix}
	    \e_k\\
        \mathbold{\f}_k
	\end{bmatrix} := 
 \begin{bmatrix}
		\nabla_\x \mathcal{L}_k(\x_k, \w_k) \\ \nabla_\w \mathcal{L}_k(\x_k, \w_k)
	\end{bmatrix} =
	\begin{bmatrix} \0 \\ \0 \end{bmatrix},
\end{equation}
where, in our case,
\begin{equation}\label{eq:gradLag}
    \begin{aligned}
        \nabla_\x \mathcal{L}_k(\x_k, \w_k) &= \A \x_k + \bv_k + \G^\top \w_k,\\
        \nabla_\w \mathcal{L}_k(\x_k, \w_k) &= \G \x_k - \h_k.
    \end{aligned}
\end{equation}
A commonly employed algorithm for solving a convex problem with linear equality constraints is based on the gradient \emph{descent-ascent} of the \emph{primal-dual} problem \cite{bertsekas_constrained_2014}. The natural extension of this \emph{primal-dual} algorithm to the online setting yields the following iterations~\cite{bernstein_online_2019,cao_online_2019}
\begin{equation}\label{eq:primalDual}
\begin{aligned}
    \x_{k+1} &= \x_k - \alpha\nabla_{\x_k} \mathcal{L}_{k}(\x_{k}, \w_{k}),\\
    \w_{k+1} &= \w_k + \beta\nabla_{\w_k} \mathcal{L}_{k}(\x_{k}, \w_{k}),
\end{aligned}
\end{equation}
where $k \in \mathbb{N}$, and $\alpha$ and $\beta$ are two positive parameters. 
Under Assumption~\ref{as:modelsbh}, we can establish that the asymptotic tracking error $\norm{\x_k-\x_k^*}$ is bounded. However, this asymptotic tracking error is not zero in general \cite{bernstein_online_2019}.

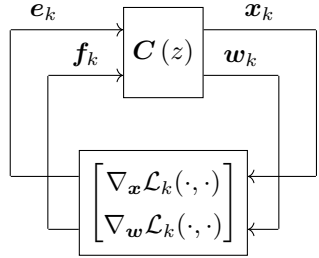
\begin{figure}[!ht]
\centering
\begin{tikzpicture}
	\node[draw,minimum height = 1.2cm] (C) at (1,0) {$\Cm\left(z\right)
		$};
	
	\node[draw=black] (grad) at (1,-2) {$\begin{bmatrix}
			\nabla_\x \mathcal{L}_k(\cdot, \cdot) \\ \nabla_\w \mathcal{L}_k(\cdot, \cdot)
		\end{bmatrix}$};
	
	\node[dot] (d1) at (-1.5,0.25) {};
    \node[dot] (d1a) at ($(C.north west)-0.25*(C.north west)+0.25*(C.south west)-(1.5,0)$) {};
    \node[dot] (d1b) at ($(C.north west)-0.75*(C.north west)+0.75*(C.south west)-(1,0)$) {};

    \node[dot] (d3a) at ($(C.north east)-0.25*(C.north east)+0.25*(C.south east)+(1.5,0)$) {};
    \node[dot] (d3b) at ($(C.north east)-0.75*(C.north east)+0.75*(C.south east)+(1,0)$) {};

    \node[dot] (d2d) at ($(grad.north west)-0.25*(grad.north west)+0.25*(grad.south west)-(-1.5,0)$) {};
    \node[dot] (d2c) at ($(grad.north west)-0.75*(grad.north west)+0.75*(grad.south west)-(0.75,0)$) {};

    \node[dot] (d4d) at ($(grad.north east)-0.25*(grad.north east)+0.25*(grad.south east)+(1.5,0)$) {};
    \node[dot] (d4c) at ($(grad.north east)-0.75*(grad.north east)+0.75*(grad.south east)+(0.75,0)$) {};

    \node[dot] (d2b) at (d1b |- d2c) {};
    \node[dot] (d2a) at (d1a |- d2d) {};
    \node[dot] (d4b) at (d3b |- d4c) {};
    \node[dot] (d4a) at (d3a |- d4d) {};
	\node[dot] (d2) at (3.25,0) {};
	\node[dot] (d3) at (3.25,-1.5) {};
	\node[dot] (d4) at (-1.5,-1.5) {};
	
    
	\draw (d1a) edge[->] node[above left] {$\mathbold{e}_k$} ($(C.north west)-0.25*(C.north west)+0.25*(C.south west)$);

    \draw (d3a) edge[-] node[above] {$\mathbold{x}_k$} ($(C.north east)-0.25*(C.north east)+0.25*(C.south east)$);
 
    \draw (d1b) edge[->] node[above] {$\mathbold{\f}_k$}  ($(C.north west)-0.75*(C.north west)+0.75*(C.south west)$);

    \path (d3b) edge[-] node[above] {$\mathbold{w}_k$} ($(C.north east)-0.75*(C.north east)+0.75*(C.south east)$);
    \path (d3a) edge (d4a)  (d4a) edge[->] ($(grad.north east)-0.25*(grad.north east)+0.25*(grad.south east)$);
    
    \path (d3b) edge (d4b)  (d4b) edge[->] ($(grad.north east)-0.75*(grad.north east)+0.75*(grad.south east)$);
    
    \path (d1b) edge (d2b)  (d2b) edge ($(grad.north west)-0.75*(grad.north west)+0.75*(grad.south west)$);
        \path (d1a) edge (d2a)  (d2a) edge ($(grad.north west)-0.25*(grad.north west)+0.25*(grad.south west)$);


		
\end{tikzpicture}
\caption{The control scheme for the solution of~\eqref{eq:quadratic-online-optimization}.}
\label{fig:block-diagram}
\end{figure}

On the other hand resorting to the control scheme illustrated in Figure~\ref{fig:block-diagram},  we can obtain an algorithm able to find all the points where \eqref{eq:error_zero} is satisfied.
The objective, therefore, is to design the transfer matrix, 
\begin{equation}\label{eq:Cmatr}
    \Cm(z) =  \begin{bmatrix}
        \Cm_{11}(z) & \Cm_{12}(z)\\
        \Cm_{21}(z) & \Cm_{22}(z)
    \end{bmatrix} \in \R^{(n+p) \times (n+p)}(z),
\end{equation} 
of the controller able to drive the signals $\e_k$ and $\f_k$ to zero.
\begin{rem}
A common alternative to the primal-dual algorithm in \eqref{eq:primalDual} is the dual-ascent algorithm where, in the update of $\w_{k+1}$, $\x_{k+1}$ is used in place of $\x_k$. In general, the performances of the dual-ascent algorithm are similar to the ones of the primal-dual algorithm. For this reason, in this paper, we decided to compare our novel strategy only with the primal-dual algorithm.
\end{rem}


\subsection{Internal Model}
Let us define the following $\mathcal{Z}$-transforms $\xh \left(z\right) \coloneqq \mathcal{Z}[\x_k]$, $\wh \left(z\right) \coloneqq \mathcal{Z}[\w_k]$, $\eh \left(z\right) \coloneqq \mathcal{Z}[\e_k]$, $\fh \left(z\right) \coloneqq \mathcal{Z}[\f_k]$.
In Figure~\ref{fig:block-diagram}, we have labeled the controller block as $ \Cm(z)$, while the input signals are denoted as time signals, following a common approach in control theory. This notation indicates that the controller represents a linear time-invariant (LTI) system acting on the input signals $\e_k$ and $\w_k$. As an LTI system, the output signals are derived from the linear combination of the convolutions of the system's impulse responses with the input signals. In terms of $\mathcal{Z}$-transforms, these relationships correspond to the following matrix multiplication.
\begin{equation}\label{eq:state}
	 \begin{bmatrix}
	    \xh \left(z\right)\\
        \wh\left(z\right)
	\end{bmatrix}= \Cm\left(z\right) \begin{bmatrix}
	    \eh \left(z\right)\\
        \fh \left(z\right)
	\end{bmatrix}.
\end{equation}
From Equation~\eqref{eq:gradLag} we also have that
\begin{equation}\label{eq:derivationCtrl}
\begin{aligned}
\begin{bmatrix}
	    \eh \left(z\right)\\
        \fh\left(z\right)
	\end{bmatrix} &= \underbrace{\begin{bmatrix}
    \A & \G^\top \\
    \G & \0
\end{bmatrix}}_{\coloneqq \mathcal{A}}\begin{bmatrix}
	    \xh \left(z\right)\\
        \wh \left(z\right)
	\end{bmatrix} + \begin{bmatrix}
	    \bh \left(z\right)\\
        -\hh\left(z\right)
	\end{bmatrix}\\
&= \mathcal{A}\Cm\left(z\right) \begin{bmatrix}
	    \eh \left(z\right)\\
        \fh\left(z\right)
	\end{bmatrix} + \begin{bmatrix}
	    \bh \left(z\right)\\
        -\hh\left(z\right)
	\end{bmatrix},
\end{aligned}
\end{equation}
where $\mathcal{A} \in \R^{(n+p) \times (n+p)}$. In this way we obtain\footnote{Where the matrix is invertible for values of $z$ different from the poles of the transfer function.}
\begin{equation}\label{eq:closeLoopSystem}
    \begin{bmatrix}
	    \eh \left(z\right)\\
        \fh\left(z\right)
	\end{bmatrix} = \left[\I - \mathcal{A}\Cm\left(z\right)\right]^{-1}\begin{bmatrix}
	    \bh \left(z\right)\\
        -\hh\left(z\right)
	\end{bmatrix}.
\end{equation}
We remark that $\mathcal{A}$ is a so-called \textit{saddle matrix} \cite{benzi_numerical_2005,benzi_eigenvalues_2006}, and in the following we will leverage results for this class of matrices to design a suitable controller.
Observe that from Assumption \ref{as:modelsbh} we know that $\bh \left(z\right),\hh\left(z\right)$ are rational and hence $\eh \left(z\right),\fh\left(z\right)$ are rational as well. Therefore, if we prove that their poles are stable, namely inside the unit circle, than we can argue that the signals $\e_k, \f_k$ converge to zero, namely
\begin{equation*}
    \begin{bmatrix}
	    \e_k\\
        \f_k
	\end{bmatrix} \xrightarrow[k \to  \infty]{}\begin{bmatrix} \0 \\ \0 \end{bmatrix}.
\end{equation*} 
By applying the Internal Model Principle~\cite{byrnes_output_1997}, the first step is to choose the controller $\Cm(z)$ able to cancel out the poles of $\bh \left(z\right),\hh\left(z\right)$, that is, the zeros of $p(z)$ introduced in Assumption \ref{as:modelsbh}. Since we assumed that all the zeros of $p(z)$ are marginally stable, then to this aim it is enough to choose 
$\Cm(z) = \frac{\Cm_\mathrm{N}(z)}{p\left(z\right)},$
where $\Cm_\mathrm{N}(z)\in \R^{(n+p) \times (n+p)}[z]$ is a polynomial matrix. Indeed, this choice yields
\begin{equation}\label{eq:error-z-transform}
	\begin{bmatrix}
	    \eh \left(z\right)\\
        \fh\left(z\right)
	\end{bmatrix} = \left[ p\left(z\right) \I - \mathcal{A} \Cm_\mathrm{N}(z) \right]^{-1} \begin{bmatrix} \bv_\mathrm{N}\left(z\right)\\ -\h_\mathrm{N}(z) \end{bmatrix}.
\end{equation}
Then the poles of $\eh \left(z\right),\fh\left(z\right)$ coincide with the poles $\left[ p\left(z\right) \I - \mathcal{A} \Cm_\mathrm{N}(z) \right]^{-1}$ and so the goal is to determine $\Cm_\mathrm{N}(z)$ such that all the poles of $\left[ p\left(z\right) \I - \mathcal{A} \Cm_\mathrm{N}(z) \right]^{-1}$ are stable.



\subsection{Stabilizing controller}

Observe that the matrix $\mathcal{A}$ is in general indefinite,
which implies that the control design approach proposed in \cite{bastianello_internal_2022} cannot be directly applied in our context.
So we have to resort to a controller $\Cm_\mathrm{N}(z)$ with a more general structure. Until now, we have considered a controller without a particular structure.
Even though imposing a structure can limit the performance of the algorithms that can be obtained, the synthesis of general controllers can be complex and intractable. For that reason, we introduce a more structured version of the controller in which we impose that
\begin{equation}\label{eq:alphaCtr}
    \Cm_\mathrm{N}\left(z\right) = c\left(z\right)\begin{bmatrix}  \I & \0 \\ \0 & -\tau\I \end{bmatrix},
\end{equation}
where $c(z) = \sum_{i = 0}^{m-1} c_i z^i \in \R[z]$ is a scalar polynomial of degree $m-1$ able to ensure strict properness of the controller transfer matrix, and $\tau$ is a positive parameter. Through this strategy we will see that it is possible to transform the design from a matrix to a scalar problem. 
With this choice of the controller we have that
\begin{equation}\label{eq:inverMatr}
\left[ p\left(z\right) \I- \mathcal{A} \Cm_\mathrm{N}(z) \right]^{-1}=
    \left[ p\left(z\right) \I -  c\left(z\right)\widetilde{\mathcal{A}}\left(\tau\right)  \right]^{-1},
\end{equation}
where 
\begin{equation}\label{eq:AcalHat}
    \widetilde{\mathcal{A}}\left(\tau\right) = \mathcal{A}\begin{bmatrix}  \I & \0 \\ \0 & -\tau\I \end{bmatrix} = \begin{bmatrix}
    \A & -\tau\G^\top \\
    \G & \0
\end{bmatrix}.
\end{equation}
We will see that the application of robust control techniques, needed in the proposed design method, is simpler if we impose that the matrix $\widetilde{\mathcal{A}}\left(\tau\right)$ has real and positive eigenvalues. The control structure described in \eqref{eq:alphaCtr} is able to guarantee this condition.

\begin{rem} One might ask whether the matrix structure of Equation~\eqref{eq:alphaCtr} can be leveraged to enhance the performance of the controller. A straightforward approach to achieving this involves introducing a new scalar parameter into the control design, resulting in the following controller:
\begin{equation*}
\Cm_\mathrm{N}\left(z\right) = c\left(z\right)\begin{bmatrix} \alpha\I & \0 \\ \0 & -\beta\I \end{bmatrix},
\end{equation*}
where $\alpha, \beta > 0$. Nevertheless, this controller can be transformed to
\begin{equation*}
    \Cm_\mathrm{N}\left(z\right) = \underbrace{c\left(z\right)\alpha}_{\coloneqq \tilde{c}\left(z\right)}\begin{bmatrix}  \I & \0  \\ \0 & -\frac{\beta}{\alpha}\I\end{bmatrix}.
\end{equation*}
Given that $c\left(z\right)$ is a polynomial that must be designed, the additional degree of freedom introduced by $\alpha$ can be integrated into this polynomial. Consequently, the problem can be addressed in the context of Equation~\eqref{eq:alphaCtr} without loss of generality.
\end{rem}
With respect to the analysis done in \cite{bastianello_internal_2022}, the design the polynomial $c(z)$ able to make all the poles of \eqref{eq:inverMatr} stable has the additional issue that matrix $\widetilde{\mathcal{A}}\left(\tau\right)$ is not always diagonalizable. Consequently, we will approach the problem using a slightly different strategy.

Consider the Schur decomposition of $\widetilde{\mathcal{A}}\left(\tau\right)$, given by $\widetilde{\mathcal{A}}\left(\tau\right) = \Q \U \Q^{\top}$, where $\Q$ is a unitary matrix and $\U$ is an upper triangular matrix. By substituting this decomposition in \eqref{eq:inverMatr} we obtain
\begin{multline}\label{eq:inverMatrSchur}
	\left[ p\left(z\right) \I -  c\left(z\right)\widetilde{\mathcal{A}}\left(\tau\right)  \right]^{-1} =\\ 
 \Q\left[ p\left(z\right) \I - c\left(z\right)\U \right]^{-1}\Q^{\top} .
\end{multline}
Hence the poles of $\left[ p\left(z\right) \I - \mathcal{A} \Cm_\mathrm{N}(z) \right]^{-1}$ coincide with the poles of $\left[ p\left(z\right) \I - c\left(z\right)\U \right]^{-1}$. Since the $\left(i,j\right)$ element of this matrix has the following form
\begin{multline}
    \label{eq:inverseMatr}
    \left(\left[ p\left(z\right) \I - c\left(z\right)\U \right]^{-1}\right)_{ij} = \\\frac{\left(\operatorname{adj}\left[ p\left(z\right) \I - c\left(z\right)\U \right]\right)_{ij}}{\det\left[ p\left(z\right) \I - c\left(z\right)\U \right]}.
\end{multline}
then the poles of this matrix are surely zeros of the polynomial $\det\left[ p\left(z\right) \I - c\left(z\right)\U \right]$. The triangular form of $\U$ implies that the set of these zeros coincides with the union of the zeros of the polynomials 
\begin{equation}\label{eq:nPoly}
p\left(z\right)-c\left(z\right)\lambda_i
\end{equation}
where $\lambda_i$, $i=1,\ldots n$, are all the eigenvalues of $\widetilde{\mathcal{A}}\left(\tau\right)$.
As emphasized in Remark~\ref{rem:robust}, our knowledge of the eigenvalues of $\A$, and consequently also of $\widetilde{\mathcal{A}}\left(\tau\right)$, is incomplete. Therefore, we need to employ a robust control technique capable to obtain stabilizing controllers of a uncertain systems. In this regard, we rely on the findings of \cite{de_oliveira_new_1999}, which employed an LMI-based approach to address this issue.

To devise a robust controller for this problem, it is essential to obtain an estimate of the eigenvalues of the matrix $\widetilde{\mathcal{A}}\left(\tau\right)$. These eigenvalues are contingent upon the specific matrices $\A$ and $\G$, as well as the parameter $\tau$. The following lemma then provides bounds for the eigenvalues of $\widetilde{\mathcal{A}}\left(\tau\right)$, and is proved by exploiting the fact that $\widetilde{\mathcal{A}}\left(\tau\right)$ is a post-conditioned saddle matrix \cite{benzi_numerical_2005,benzi_eigenvalues_2006}.

\begin{lem}[Eigenvalues of $\widetilde{\mathcal{A}}\left(\tau\right)$]\label{lem:eigenvalues_A_hat} 
Given $\widetilde{\mathcal{A}}\left(\tau\right)$ as in \eqref{eq:AcalHat}, if $\tau$ is such that
\begin{equation}\label{eq:optChoiceineq}
    0<\tau \leq \tau^* \coloneqq \frac{\lmin\left(\A\right)}{4\lmax\left(\G\A^{-1}\G^{\top}\right)}.
\end{equation} 
then the eigenvalues of $\widetilde{\mathcal{A}}\left(\tau\right)$ are real and belong to the following interval
\begin{equation}\label{eq:interval}
\left[ \tau\lmin\left(\G\A^{-1}\G^{\top}\right), \lmax\left(\A\right)\right]
\end{equation}
\end{lem}
\begin{pf}
As proved in \cite[Corollary 2.7]{benzi_eigenvalues_2006} if $\lmin\left( \A\right)\geq 4\tau\lmax\left(\G\A^{-1}\G^{\top}\right)$, then $\widetilde{\mathcal{A}}\left(\tau\right)$ has all real eigenvalues, which implies inequality~\eqref{eq:optChoiceineq} and this implies that $\widetilde{\mathcal{A}}\left(\tau\right)$ has all real eigenvalues.
Moreover, thanks to the fact that $\G$ is full row rank and using \cite[Proposition 2.12]{benzi_eigenvalues_2006} and \cite[Theorem 2.1]{shen_eigenvalue_2010} we know that the eigenvalues of $\widetilde{\mathcal{A}}\left(\tau\right)$ belongs to the interval 
$$\left[ \min \left\lbrace\lmin\left(\A\right),\tau\lmin\left(\G\A^{-1}\G^{\top}\right)\right\rbrace, \lmax\left(\A\right)\right]$$
Observe finally that, using \eqref{eq:optChoiceineq} we can argue that\footnote{Recalling that $\kappa$ denotes the condition number.}
\begin{equation}\label{eq:ineqInterval1}
\begin{aligned}
        \tau\lmin\left(\G\A^{-1}\G^{\top}\right)&\le\frac{\lmin\left(\A\right)}{4\lmax\left(\G\A^{-1}\G^{\top}\right)}\lmin\left(\G\A^{-1}\G^{\top}\right) \\
        &= \frac{\lmin\left(\A\right)}{4\kappa\left(\G\A^{-1}\G^{\top}\right)}\leq \lmin\left(\A\right),
\end{aligned}
\end{equation}
which yields \eqref{eq:interval}. 
\qed
\end{pf}
    
    Based on Lemma~\ref{lem:eigenvalues_A_hat}, we can determine $\tau$ such that the eigenvalues of $\widetilde{\mathcal{A}}\left(\tau\right)$ are real and fall in the interval specified in Equation~\eqref{eq:interval}. To stabilize the polynomials given by Equation~\eqref{eq:nPoly}, it is sufficient to select an appropriate $c\left(z\right)$ that stabilizes the polynomials
\begin{equation}\label{eq:staPol}
    p\left(z\right)-c\left(z\right)\lambda,
\end{equation}
for all $\lambda \in \left[\lmin\left(\widetilde{\mathcal{A}}\left(\tau\right)\right),\lmax\left(\widetilde{\mathcal{A}}\left(\tau\right)\right)\right]$. But requiring
the stability of the polynomials in \eqref{eq:staPol} coincides with requiring the stability of the associated companion matrices $\F_c\left(\lambda\right)\coloneqq \F + \lambda \mathbold{C} \K$, where

\begin{equation}\label{eq:compMatr}
\begin{aligned}
      \F &= \begin{bNiceMatrix}
0      &   1     & 0 &   \Cdots  & 0 \\
      & \Ddots  &  \Ddots    & \Ddots   & \Vdots \\
\Vdots &   &        &      & 0 \\
0      & \Cdots  &       & 0     & 1 \\
-p_0   &         & \Cdots &      & -p_{m-1}    
\end{bNiceMatrix} \quad
      \mathbold{C}&=\begin{bNiceMatrix}
        0 \\
          \Vdots\\
         \\
         \\
         0\\
         1
    \end{bNiceMatrix} 
      \\
      \K &= \begin{bNiceMatrix}[columns-width=5.6mm] c_0 & & \Cdots &  &c_{m-1} \end{bNiceMatrix} & .
\end{aligned}
\end{equation}

Realization~\eqref{eq:compMatr} represents a state-space realization of the transfer function $c\left(z\right)/p\left(z\right)$ that constitutes the controller $\Cm(z)$ in Figure \ref{fig:block-diagram}, since we have that
\begin{equation*}
    \Cm\left(z\right) = \frac{c\left(z\right)}{p\left(z\right)}\begin{bmatrix}  \I & \0 \\ \0 & -\tau\I \end{bmatrix},
\end{equation*}
In this way the zeros of the polynomial in \eqref{eq:nPoly} coincide with the poles of the closed loop transfer matrix.
By denoting $\ubar{l} \coloneqq \lmin\left(\widetilde{\mathcal{A}}\left(\tau\right)\right)$ and $\bar{l} \coloneqq \lmax\left(\widetilde{\mathcal{A}}\left(\tau\right)\right)$, we can express $\lambda$ as a convex combination of these two extreme values using the equation:
\begin{equation*}
    \lambda = \alpha\left(\lambda\right)\ubar{l} + \left(1-\alpha\left(\lambda\right)\right)\bar{l}, \quad \alpha\left(\lambda\right) = \frac{\bar{l}-\lambda}{\bar{l}-\ubar{l}}
\end{equation*}
With this expression, we can apply the following result \cite[Theorem 3]{de_oliveira_new_1999}), and subsequently derive the polynomial $c\left(z\right)$ by solving the two LMIs of size $m$ in the following lemma.

\begin{lem}[Stabilizing controller design]\label{lem:stabC}
The eigenvalues of $\F_c\left(\lambda\right) $ are strictly inside the unit circle for any $\lambda \in \left[\,\ubar{l},\bar{l}\,\right]$ if there exist symmetric matrices $\underline{\P}, \overline{\P} \succ 0$, a square matrix $\Q\in \R^{m\times m }$, and a row vector $\mathbold{R}\in\mathbb{R}^{1 \times m}$ such that
\begin{align*}
    \begin{bmatrix}
        \underline{\P} & \F\Q + \ubar{l}\mathbold{C}\mathbold{R}\\
        \Q^{\top}\F^{\top}+\ubar{l}\mathbold{R}^{\top}\mathbold{C}^{\top} & \Q + \Q^{\top} - \underline{\P}
    \end{bmatrix} &\succ 0,\\
    \begin{bmatrix}
        \overline{\P} & \F\Q + \bar{l}\mathbold{C}\mathbold{R}\\
        \Q^{\top}\F^{\top}+\bar{l}\mathbold{R}^{\top}\mathbold{C}^{\top} & \Q + \Q^{\top} - \overline{\P}
    \end{bmatrix} &\succ 0.
\end{align*}
A stabilizing controller is then $\K = \mathbold{R}\Q^{-1}$.
\end{lem} These LMIs depend on the minimum and maximum eigenvalues of $\widetilde{\mathcal{A}}\left(\tau\right)$. But, using Lemma \ref{lem:eigenvalues_A_hat}, we can obtain the polynomial $c\left(z\right)$ by solving the two LMIs on the extremes of the interval in \eqref{eq:interval}.

The feasibility of aforementioned robust control design procedure, based on \cite[Theorem 3]{de_oliveira_new_1999} depends on the ratio between the extremes of the interval. Precisely the possibility of obtaining a solution increases as much as this ratio get closer to one.
Hence the best choice of the parameter $\tau$ is the maximum value ensuring that the eigenvalues of $\widetilde{\mathcal{A}}\left(\tau\right)$ are real. 
From Lemma \ref{lem:eigenvalues_A_hat} this is $\tau=\tau^*$ and in this case we can argue that the eigenvalues of $\widetilde{\mathcal{A}}\left(\tau^*\right)$ belong to the interval
$$\left[\frac{\lmin\left(\A\right)\lmin\left(\G\A^{-1}\G^{\top}\right)}{4\lmax\left(\G\A^{-1}\G^{\top}\right)},\lmax\left(\A\right)\right],$$
and that the ratio between the extremes of this interval is
$$4\kappa\left(\A\right)\kappa\left(\G\A^{-1}\G^{\top}\right).$$
In the case where only the bounds given in Equations \eqref{eq:estA} and \eqref{eq:estGAGT} are available, we select $\tau=\tau^*_{\mathrm{est}} := \lmin/4\ubar{\mu}$. Consequently, the eigenvalues of $\widetilde{\mathcal{A}}\left(\tau^*_{\mathrm{est}}\right)$ fall within the interval $
    \left[\lmin\ubar{\mu}/ \left(4\bar{\mu}\right),\lmax\right]$ and the ratio between the extremes of this interval is $4\lmax\bar{\mu}/\left(\lmin\ubar{\mu}\right)$.

Finally, we remark that the computational complexity of solving these LMIs depends also on the degree of $p(z)$. However, the solution of the LMIs is performed only once, offline before the deployment of the algorithm. Thus the LMIs complexity does not affect the real-time operation.

\subsection{The online algorithm}

The algorithm involves translating the relationships among the signals $\e_k$, $\f_k$, $\x_k$, and $\w_k$ from the $\mathcal{Z}$-transform domain to the time domain. These relationships are defined by the following equations in the $\mathcal{Z}$-transform domain:
$$\xh \left(z\right)=\left(\frac{c(z)}{p(z)}\I\right)\eh \left(z\right),\qquad \wh \left(z\right)=-\tau\left(\frac{c(z)}{p(z)}\I\right)\fh \left(z\right).$$ Specifically, the algorithm's recursion can be expressed by directly deriving a realization of the transfer matrices of the previous controllers, employing a standard approach in control theory.
A state space realization of such relation is given by
\begin{subequations}\label{eq:control-based-algorithm}
\begin{empheq}[box=\widefbox]{align}
	\z_{k+1} &= \left(\F \otimes \I\right)  \z_k + \left(\Cm \otimes \I \right) \e_k\label{eq:control-based-algorithm-y} \\
    \y_{k+1} &= \left(\F \otimes \I\right)  \y_k + \left(\Cm \otimes \I \right) \f_k\label{eq:control-based-algorithm-z} 
 	\end{empheq} 
 \begin{empheq}[box=\widefbox]{align}
	\x_{k+1} &= \left( \K \otimes \I \right) \z_{k+1}\label{eq:control-based-algorithm-x}\\
	\w_{k+1} &= - \tau\left( \K \otimes \I\right) \y_{k+1} \label{eq:control-based-algorithm-w}
    \end{empheq} 
\begin{empheq}[box=\widefbox]{align}
    \e_{k+1}&=\nabla_\x \mathcal{L}_k(\x_{k+1}, \w_{k+1}) \label{eq:control-based-algorithm-e} \\
    \f_{k+1}&=\nabla_\w \mathcal{L}_k(\x_{k+1}, \w_{k+1})\label{eq:control-based-algorithm-f} 
\end{empheq} 
\end{subequations}
where the matrices $\F$ , $\Cm$ and $\K$ are defined in \eqref{eq:compMatr} and $\z_k \in \R^{nm},\y_k \in \R^{pm}$ are the state vectors of the controller. We remark that we use the Kronecker products of matrices $\F$, $\Cm$ and $\K$ with $\I \in \R^{n \times n}$ for~\eqref{eq:control-based-algorithm-y} and~\eqref{eq:control-based-algorithm-x} to account for the fact that the optimization variable $\x$ has dimension $n$. Similarly, the Kronecker products in Equations~\eqref{eq:control-based-algorithm-z} and~\eqref{eq:control-based-algorithm-w} are done with $\I \in \R^{p \times p}$ accounts for the fact that we have $p$ constraints.
Hence, equations \eqref{eq:control-based-algorithm} describe the online algorithm in the scheme of Figure~\ref{fig:block-diagram}.


\subsection{Convergence}\label{subsec:convergence}

In this section we prove the convergence of Algorithm~\eqref{eq:control-based-algorithm} solving the optimization problem described in \eqref{eq:quadratic-online-optimization}.
Convergence is established by the following proposition, which employs classical arguments based on a $\mathcal{Z}$-transforms.

\begin{prop}[Convergence of Algorithm~\ref{eq:control-based-algorithm}]\label{pr:convergence-control-algorithm}
Consider the optimization problem~\eqref{eq:quadratic-online-optimization}, with $\{ \bv_k \}_{k \in \N}$ and $\{ \h_k \}_{k \in \N}$ modeled by~\eqref{eq:linearModels}. Choose the controller~\eqref{eq:alphaCtr} with $\tau$ satisfying \eqref{eq:optChoiceineq}
and let $c(z)$ be such that $\F_{\mathrm{c}}(\lambda)$, defined in~\eqref{eq:compMatr}, is asymptotically stable for all $\lambda \in\left[ \tau\lmin\left(\G\A^{-1}\G^{\top}\right), \lmax\left(\A\right)\right]$. Let $\{ \x_k \}_{k \in \N}$ be the output of the online algorithm~\eqref{eq:control-based-algorithm}. Then it holds
$$
	\limsup_{k \to \infty} \norm{\x_k - \x_k^*} = 0.
$$
\end{prop}
\begin{pf}
For the proof we follow the reasoning of \cite[Proposition 1]{bastianello_internal_2022}. Given the control structure~\eqref{eq:alphaCtr} and considering that interval $\left[ \tau\lmin\left(\G\A^{-1}\G^{\top}\right), \lmax\left(\A\right)\right]$ includes the interval $\left[\lmin\left(\widetilde{\mathcal{A}}\left(\tau\right)\right), \lmax\left(\widetilde{\mathcal{A}}\left(\tau\right)\right)\right]$, selecting a controller $c\left(z\right)$ that stabilizes the matrix $\F_{\mathrm{c}}\left(\lambda\right)$ for $\lambda \in\left[ \tau\lmin\left(\G\A^{-1}\G^{\top}\right), \lmax\left(\A\right)\right]$ ensures the asymptotic stability of the poles of $\eh \left(z\right)$ and of $\fh\left(z\right)$, see~\eqref{eq:error-z-transform}. Consequently, the gradient of the Lagrangian, given by $\e_k$, $\f_k$, converges to zero, and this implies the thesis. \qed
\end{pf}

\begin{rem}\label{rem:genCase}
If the cost has a time-varying $\A$ or is non-quadratic\footnote{Similarly, if the constraints have $\G$ time-varying or are non-linear.}, the convergence of the tracking error to zero is not guaranteed. This is because our approach reformulates problem~\eqref{eq:quadratic-online-optimization-general} as a linear control problem, while it becomes a non-linear control problem for more general classes of problems.
To analyze the performance for these problems, one can resort to the small gain theorem, which however tends to over-estimate the tracking error.
Alternatively, one can design the online algorithm using the non-linear internal model principle.
We leave these theoretical developments for future research, and in this paper we only provide (promising) numerical evaluations of the algorithm for general problems.

\end{rem}

\section{Online Optimization with Equality and Inequality Constraints}\label{sec:ineq-constraints}

In this section, we aim to extend the previous approach to also handle linear inequality constraints.

We now consider Problem~\eqref{eq:quadratic-online-optimization-general} addressing both equality and inequality constraints simultaneously. Specifically, as outlined in Assumption~\ref{as:modelsbh} we assume partial knowledge of the time-varying terms:
\begin{assum}[Models of $\bv_k$, $\h_k$ and $\h^{\prime}_k$]\label{as:modelsbhhprime}
The sequence of vectors $\left\lbrace\bv_k \right\rbrace_{k \in \N}$,  $\left\lbrace \h_k \right\rbrace_{k \in \N}$ and $\left\lbrace \h^{\prime}_k \right\rbrace_{k \in \N}$ have rational $\mathcal{Z}$-transforms
\begin{align}
    \bh\left(z\right) \coloneqq \mathcal{Z}\left[ \bv_k \right]& = \frac{\bv_\mathrm{N} \left(z\right)}{p\left(z\right)}, \quad
    \hh\left(z\right) \coloneqq \mathcal{Z}\left[ \h_k \right] = \frac{\h_\mathrm{N}\left(z\right)}{p\left(z\right)} , \nonumber\\
    &\hh^{\prime}\left(z\right) \coloneqq\mathcal{Z}\left[ \h_k^{\prime} \right] = \frac{\h_\mathrm{N}^{\prime} \left(z\right)}{p\left(z\right)} \label{eq:linearModels}
\end{align}
where the polynomial $p\left(z\right) = z^m + \sum_{i = 0}^{m-1} p_i z^i$ is known and whose zeros are assumed to be all marginally stable\footnote{The same reasoning presented in Remark~\ref{rem:assum} also applies to Assumption~\ref{as:modelsbhhprime}.}. The numerators, $\bv_\mathrm{N}\left(z\right)$, $\h_\mathrm{N}\left(z\right)$ and $\h_\mathrm{N}^{\prime}\left(z\right)$, are instead assumed to be unknown.
\end{assum}
A commonly employed approach to address this problem is to extend the \emph{primal-dual} algorithm to the (online) \emph{projected primal-dual}~\cite{cao_online_2019}, which yields
\begin{equation}\label{eq:proPrimalDual}
\begin{aligned}
    \x_{k+1} &= \x_k - \alpha\nabla_{\x_k} \mathcal{L}^{\prime}_{k}(\x_{k}, \w_{k},\w_{k}^{\prime}),\\
    \w_{k+1} &= \w_k + \beta\nabla_{\w_k} \mathcal{L}^{\prime}_{k}(\x_{k}, \w_{k},\w_{k}^{\prime}),\\
    \w_{k+1}^{\prime} &= \proj_{\geq 0}\left(\w^{\prime}_k + \gamma\nabla_{\w^{\prime}_k} \mathcal{L}^{\prime}_{k}(\x_{k}, \w_{k},\w_{k}^{\prime})\right).
    \end{aligned}
\end{equation}
Here, $\alpha$, $\beta$, and $\gamma$ are positive parameters, $\w_k \in \mathbb{R}^{p}$ and $\w_k^\prime \in \mathbb{R}^{p^\prime}$ are the Lagrange multipliers associated with the equality and the inequality constraints and $\mathcal{L}_{k}^{\prime}$ is the time-varying Lagrangian defined as
\begin{multline*}
    \mathcal{L}_{k}^{\prime}\left(\x,\w,\w^{\prime}\right) \coloneqq\\ f_k\left(\x\right) + \w^{\top}\left(\G\x-\h_k\right) + {\w^{\prime}}^{\top}\left(\G^{\prime}\x-\h^{\prime}_k\right).
 \end{multline*}

The main difference between Algorithms~\eqref{eq:primalDual} and \eqref{eq:proPrimalDual} is the presence of a projection in the update of the dual variable associated with the inequality constraints. This projection works as a saturation to zero of the negative entries of the signal. 

\begin{figure}[!ht]
\centering
\begin{tikzpicture}
	\node[draw,minimum width
=1.5cm, minimum height = 1.8cm] (C) at (1,0) {Equations (31a) - (31f)};
	
	\node[draw=black,minimum size = 2cm] (grad) at (1,-3.5) {$\begin{bmatrix}
			\nabla_{\x} \mathcal{L}^{\prime}_{k}(\cdot, \cdot,\cdot) \\ \nabla_{\w} \mathcal{L}^{\prime}_{k}(\cdot, \cdot,\cdot)\\
   \nabla_{\w^{\prime}} \mathcal{L}^{\prime}_{k}(\cdot, \cdot,\cdot)
		\end{bmatrix}$};
 
	\node[dot] (d1) at (-1.5,0.25) {};
    \node[dot] (d1a) at ($(C.north west)-0.25*(C.north west)+0.25*(C.south west)-(1.05,0)$) {};
    \node[dot] (d1b) at ($(C.north west)-0.5*(C.north west)+0.5*(C.south west)-(0.65,0)$) {};
    \node[dot] (d1c) at ($(C.north west)-0.75*(C.north west)+0.75*(C.south west)-(0.25,0)$) {};
    
    \node[dot] (d3a) at ($(C.north east)-0.25*(C.north east)+0.25*(C.south east)+(3.0,0)$) {};
    \node[dot] (d3b) at ($(C.north east)-0.5*(C.north east)+0.5*(C.south east)+(2.8,0)$) {};
    \node[dot] (d3c) at ($(C.north east)-0.75*(C.north east)+0.75*(C.south east)+(2.6,0)$) {};

    \node (d5a) at ($(grad.north west)-0.25*(grad.north west)+0.25*(grad.south west)-(1.75,0)$) {};
    \node (d5b) at ($(grad.north west)-0.5*(grad.north west)+0.5*(grad.south west)-(1.15,0)$) {};
    \node (d5c) at ($(grad.north west)-0.75*(grad.north west)+0.75*(grad.south west)-(0.45,0)$) {};
    \node[dot] (d6b) at ($(grad.north west)-0.75*(grad.north west)+0.75*(grad.south west)-(0.5,0)$) {};

    \node[dot] (d2d) at ($(grad.north west)-0.25*(grad.north west)+0.25*(grad.south west)-(-1.5,0)$) {};
    
    \node[dot] (d2in) at ($(grad.north west)-0.75*(grad.north west)+0.75*(grad.south west)-(0.75,0)$) {};
    \node[dot] (d2in0) at ($(grad.north west)-0.5*(grad.north west)+0.5*(grad.south west)-(0.75,0)$) {};
    \node[dot] (d4d) at ($(grad.north east)-0.25*(grad.north east)+0.25*(grad.south east)+(1.5,0)$) {};
    \node[dot] (d4in1) at ($(grad.north east)-0.5*(grad.north east)+0.5*(grad.south east)+(1.5,0)$) {};
    \node[dot] (d4in0) at ($(grad.north east)-0.75*(grad.north east)+0.75*(grad.south east)+(0.75,0)$) {};

    \node[dot] (d2c) at (d1c |- d2in) {};
    \node[dot] (d2b) at (d1b |- d2in0) {};
    \node[dot] (d2a) at (d1a |- d2d) {};
    \node[dot] (d4c) at (d3c |- d4in0) {};
    \node[dot] (d4b) at (d3b |- d4in1) {};
    \node[dot] (d4a) at (d3a |- d4d) {};
	\node[dot] (d2) at (3.25,0) {};
	\node[dot] (d3) at (3.25,-1.5) {};
	\node[dot] (d4) at (-1.5,-1.5) {};
  
	\draw (d1a) edge[->] node[above left,pos=0.7,yshift=-2pt] {$\mathbold{e}_k$} ($(C.north west)-0.25*(C.north west)+0.25*(C.south west)$);

    \draw (d3a) edge[-] node[above right,pos=0.25,yshift=-2pt] {$\mathbold{x}_k$} ($(C.north east)-0.25*(C.north east)+0.25*(C.south east)$);

    \draw (d3b) edge[-] node[above right,pos = 0.25, yshift = -2pt] {$\w_k$} ($(C.north east)-0.5*(C.north east)+0.5*(C.south east)$);
 
    \draw (d1b) edge[->] node[above,pos=0.3,yshift=-2pt] {$\f_k$}  ($(C.north west)-0.5*(C.north west)+0.5*(C.south west)$);

    \draw (d1c) edge[->] node[above,pos=-0.3,yshift=-2pt] {$\f^{\prime}_k$}  ($(C.north west)-0.75*(C.north west)+0.75*(C.south west)$);

    \path (d3c) edge[-] node[above right,pos=0.25,yshift= -2pt] {$\w^{\prime}_k$} node[above left,yshift= -2pt,pos=0.75] {$\vv_k$}($(C.north east)-0.75*(C.north east)+0.75*(C.south east)$)  ;
    \path (d3a) edge (d4a)  (d4a) edge[->] ($(grad.north east)-0.25*(grad.north east)+0.25*(grad.south east)$);
    
    \path (d3b) edge (d4b)  (d4b) edge[->] ($(grad.north east)-0.5*(grad.north east)+0.5*(grad.south east)$);
    
    \path (d3c) edge (d4c)  (d4c) edge[->] ($(grad.north east)-0.75*(grad.north east)+0.75*(grad.south east)$);
    
    \path (d1c) edge (d2c)  (d2c) edge ($(grad.north west)-0.75*(grad.north west)+0.75*(grad.south west)$);
    \path (d1b) edge (d2b)  (d2b) edge ($(grad.north west)-0.5*(grad.north west)+0.5*(grad.south west)$);
        \path (d1a) edge (d2a)  (d2a) edge ($(grad.north west)-0.25*(grad.north west)+0.25*(grad.south west)$);
    
    \node[saturation block,minimum width=0.7cm,minimum height=0.7cm,,draw] (sat) at ($(C.north east) - 0.75*(C.north east) + 0.75*(C.south east)+(1.3,0)$) {};
    \node[dot] (satTmp) at ($(sat) - (0,0.75)$) {};
    \node (d6b) at (sat |- satTmp) {};
    \node[dot] (tmpd6b) at ($(d6b) - (0.75,0)$) {};
    \draw[fill=white](d6b) circle (0.05) edge [<-] (tmpd6b) node [below left] {$-$};
    
    \node[dot] (tmpd7b) at ($(d6b) + (0.75,0)$) {};
    \draw(d6b) edge [<-] (tmpd7b) node [below right] {$+$};
    \path (tmpd7b) edge (tmpd7b|-d3c);
    \path (tmpd6b) edge (tmpd6b|-d3c);
    
    \node[dot] (tmpd8b) at ($(d6b) - (0,0.6)$) {};
    
    \path (d6b) edge (tmpd8b);
    \node (tmpd9b) at ( d1c|- tmpd8b) {};
    \draw[fill=white]( tmpd9b) circle (0.05) edge [<-] (tmpd8b) node [below right] {$+$};
    \node[dot] (tmp) at ($(d5a) - (0,1.75)$) {};
	\node[dot] (tmp) at ($(d5c) - (0,0.75)$) {};

    \node[dot] (tmp) at ($(d5b) - (0,1.25)$) {};
	\node[dot] (tmpRho) at (tmpd8b-|grad) {};
    \node[fill=white,draw=black,minimum size =0.5cm] (rho) at (tmpRho) {$\rho$};
		
\end{tikzpicture}
\caption{The control scheme designed to solve~\eqref{eq:quadratic-online-optimization-general}.}
\label{fig:block-diagram-anti wind-up}
\end{figure}
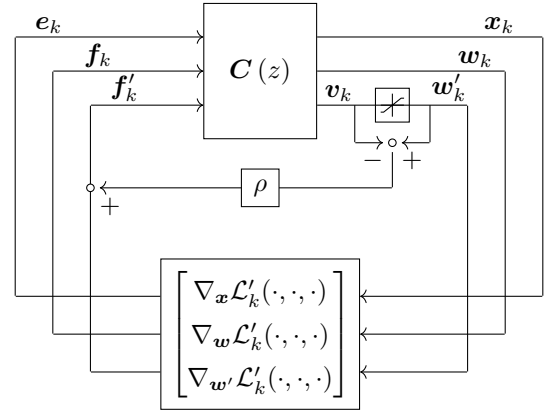

The first attempt in the design of an online optimization algorithm based on the Internal Model principle is simply applying what we did in the case in which we have only equality constraints. Precisely, we determine
the controller using the internal model $p\left(z\right)$ as done in the previous section, considering inequality constraints as they were further equality constraints, and we incorporate in the resulting algorithm the saturation as done in \eqref{eq:proPrimalDual}. 
However, if we do this, we observe in the experiments the output of the proposed algorithm achieves very poor tracking, worse than the unstructured online primal-dual algorithm. Indeed, it is well-known that a controller based on the internal model with marginally stable poles, coupled with a saturation, gives rise to significant transient phenomena known as wind-up. In order to mitigate this negative effect of the saturation, an anti wind-up mechanism, called back-calculation \cite{astrom_antiwindup_1989}, is added to the algorithm\footnote{Other approaches to mitigate the wind-up effect exist, such as those discussed in~\cite{Tarbouriech_antiwindup_2009}; however, these methods are not covered in this work.}. 
The resulting algorithm is illustrated in Figure~\ref{fig:block-diagram-anti wind-up} and described by the following equations~\eqref{eq:control-based-algorithma-anti wind-up}
\begin{subequations}\label{eq:control-based-algorithma-anti wind-up}
\begin{empheq}[box=\widefbox]{align}
	\z_{k+1} &= \left(\F \otimes \I\right)  \z_k + \left(\Cm \otimes \I \right) \e_k\label{eq:control-based-algorithm-y1} \\
    \y_{k+1} &= \left(\F \otimes \I\right)  \y_k + \left(\Cm \otimes \I \right) \f_k\label{eq:control-based-algorithm-z1} \\
    \y_{k+1}^{\prime} &= \left(\F \otimes \I\right)  \y_k^{\prime} + \left(\Cm \otimes \I \right) \f_k^{\prime}\label{eq:control-based-algorithm-zp} 
    \end{empheq} 
\begin{empheq}[box=\widefbox]{align}
    \x_{k+1} &= \left( \K \otimes \I \right) \z_{k+1}\\\
	\w_{k+1} &= - \tau\left( \K \otimes \I \right) \y_{k+1}\\
    \v_{k+1} &= - \tau\left( \K \otimes \I \right) \y_{k+1}^{\prime}
\end{empheq} 
\begin{empheq}[box=\fbox]{align}
 \w_{k+1}^{\prime} &= \proj_{\geq 0}\left(\v_{k+1}\right)\label{eq:satanti wind-up}\\
 \e_{k+1}&= \nabla_\x \mathcal{L}^{\prime}_{k+1}(\x_{k+1}, \w_{k+1},\w_{k+1}^{\prime})\\
 \f_{k+1} &= \nabla_{\w} \mathcal{L}^{\prime}_{k+1}(\x_{k+1}, \w_{k+1},\w_{k+1}^{\prime})\\
  \f_{k+1}^{\prime} &= \nabla_{\w^{\prime}} \mathcal{L}^{\prime}_{k+1}(\x_{k+1}, \w_{k+1},\w_{k+1}^{\prime}) +\\ &\qquad+\rho\underbrace{\left(\w^{\prime}_{k+1}-\v_{k+1}\right)}_{\mathrm{anti wind-up}}
.\label{eq:errorDualanti wind-up}
\end{empheq}
\end{subequations}

where $\rho$ is a positive tunable anti wind-up parameter, the matrices $\F$, $\Cm$ and $\K$ are defined in \eqref{eq:compMatr} and $\z_k \in \R^{nm},\y_k \in \R^{pm},\y_k ^{\prime}\in \R^{p^{\prime}m}$ are the state vectors of the controller.

Analyzing the tracking error of the algorithm in the case of inequality constraints, both with and without anti-windup, is not straightforward. In particular, transient phenomena occur whenever inequality constraints switch from active to inactive, precluding exact tracking without additional (restrictive) assumptions. We illustrate and discuss this behavior in section~\ref{subsec:sim-inequality-constraints}.


\section{Numerical Experiments}\label{sec:simulations}
In this section we compare the performance of the different proposed algorithms when applied to the Problem~\eqref{eq:quadratic-online-optimization} and to the Problem~\eqref{eq:quadratic-online-optimization-general}. In addition to addressing Problem~\eqref{eq:quadratic-online-optimization}, we also provide a numerical example to illustrate that the algorithm performs well even when the matrices $\A$ and $\G$ are time-varying and in more general non-quadratic cases.

\subsection{Equality constraints}\label{subsec:sim-equality-constraints}

First we consider Problem~\eqref{eq:quadratic-online-optimization} with $\x  \in \mathbb{R}^n$, where $n = 10$. The matrix $\A$ is defined as $\A = \V\Lm\V^{\top}$, where $\V$ is a randomly generated orthogonal matrix, and $\Lm$ is a diagonal matrix with elements in the range $[1,10]$. On the other hand, $\G \in \mathbb{R}^{p\times n}$ is a randomly generated matrix with orthogonal rows.

\begin{figure}[!ht]
\centering
\includegraphics[scale = 0.5]{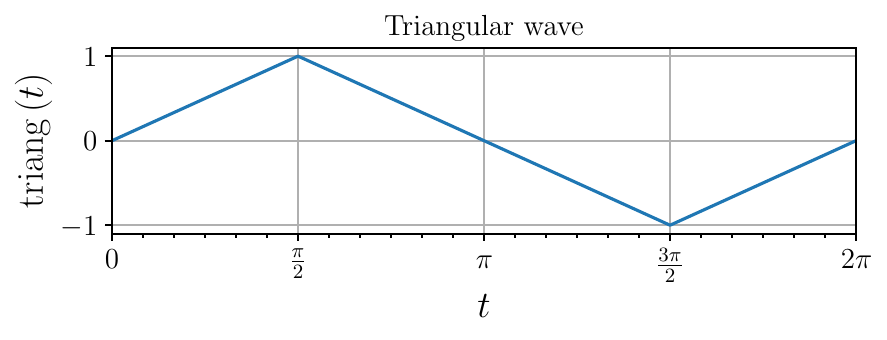}
\caption{The graph of the $\mathrm{triang}\left(t\right)\coloneqq4 \left|\frac{t}{2\pi}-\lfloor \frac{t}{2\pi} + 3/4\rfloor+1/4\right|-1$,  $t\in \R$.}
\label{fig:triang}
\end{figure}

\begin{figure}[!ht]
\centering
\includegraphics[scale = 0.5]{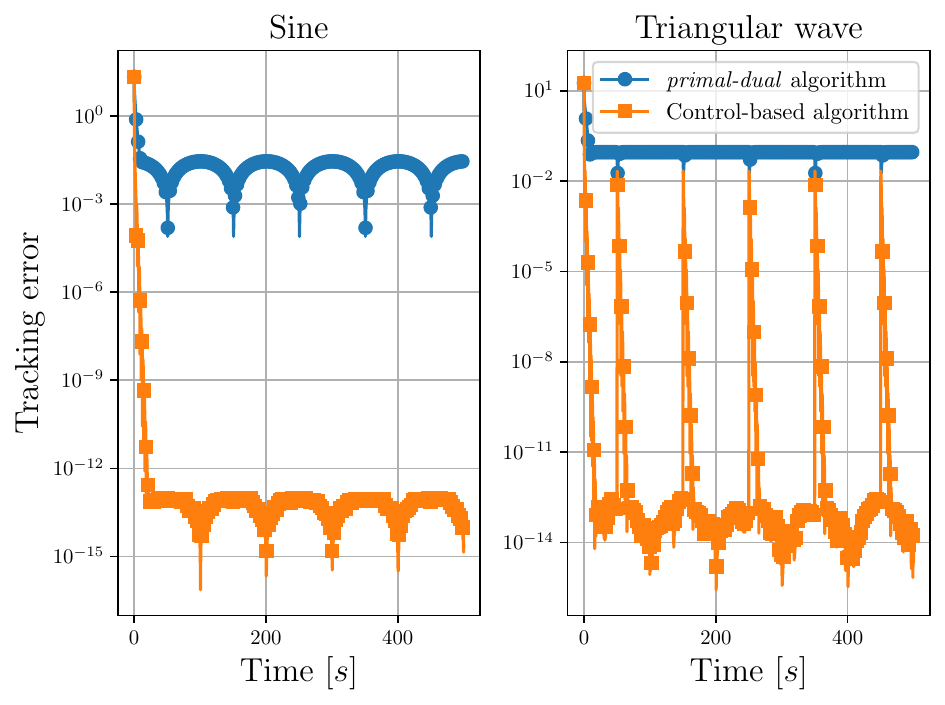}
\caption{Comparison between the \emph{primal-dual} and the control-based algorithm in a semilogarithmic plot in case of sinusoidal or triangular wave signals.}
\label{fig:comp-eq}
\end{figure}

For the linear terms $\left\lbrace\bv_k \right\rbrace_{k \in \N}$ and $\left\lbrace \h_k \right\rbrace_{k \in \N}$, we used the following models:
\begin{enumerate}
    \item Triangular wave (Figure~\ref{fig:triang}): $\bv_k = \mathrm{triang}\left(\omega k\right)\1$, $\h_k = \mathrm{triang}\left(\omega k\right)\1$ with $\omega=10^{-4}\pi$;
    \item Sine: $\bv_k = \sin\left(\omega k\right)\1$,  $\h_k = \sin\left(\omega k\right)\1$ with $\omega=10^{-4}\pi$. \label{enum:case2-eq}
\end{enumerate}
In Figure~\ref{fig:comp-eq}, we compare our algorithm, which is based on a control approach with $p\left(z\right) = z^2 - 2\cos\left(\omega\right)z + 1$ for the sine and $p\left(z\right) = \left(z - 1\right)^2$ for the triangular wave, to the standard \emph{primal-dual} \cite{cao_online_2019} algorithm given by~\eqref{eq:primalDual}. From a control perspective, the classical \emph{primal-dual} algorithm can be interpreted as a controller with $p\left(z\right) = \left(z - 1\right)$ as its internal model. Consequently, for the two waves, the \emph{primal-dual} algorithm fails to ensure perfect tracking because it does not satisfy the internal model principle. On the other hand, the control approach design based on the correct model is capable of guaranteeing perfect tracking. In the case of the triangular wave, we can see that the Internal Model principle is piece-wise verified and then the tracking is always guaranteed except when there is a change of the slope in the signal. In these time instants a transient is required by  the controller in order to track back the reference.

We conclude this section by presenting a simulation that addresses the scenario in which the model is ``slightly'' incorrect. Specifically, we assume that $\bv_k$ and $\h_k$ follow the same evolution as in the ``Sine'' case depicted in Fig.~\ref{fig:comp-eq}. We conduct a series of simulations to measure the asymptotic error for an algorithm designed with an incorrect oscillation frequency $\tilde{\omega}$ compared to the actual oscillation frequency $\omega$. The results are illustrated in Table~\ref{tab:tabInexModel}. The data illustrates that, although perfect tracking is not achieved, the control-based algorithm outperforms the \emph{primal-dual} algorithm even when the design frequency deviates from the actual frequency by up to $20\%$.
\begin{table}
\[
\begin{array}{cccccc}
\toprule
\frac{|\omega-\tilde{\omega}|}{\omega}= & 4\% & 8\% & 12\% & 16\% & 20\% \\
\midrule
\text{\emph{primal-dual} algorithm} & 3.73 & - & - & - & -  \\
\text{control-based algorithm}& 0.20 & 0.41 & 0.61 & 0.82 & 1.02\\
\bottomrule
\end{array}
\]
\caption{Asymptotic errors of the \emph{primal-dual} algorithm compared to the control-based algorithm when $\bv_k$ is given by model~\ref{enum:case2-eq} (the ``Sine'' model) as already described in Section~\ref{subsec:sim-equality-constraints} with $p(z) = z^2 - 2\cos(\omega)z + 1$. The control-based algorithm in~\eqref{eq:control-based-algorithm} is designed with $p(z) = \tilde{p}(z) \coloneqq z^2 - 2\cos(\tilde{\omega})z + 1$ where $\tilde{\omega} = 0.96 \omega, 0.92 \omega, 0.88 \omega, 0.84 \omega, 0.8 \omega$. The \emph{primal-dual} algorithm yields a single error value since its design is independent of the oscillation frequency of $\bv_k$.}
\label{tab:tabInexModel}
\end{table}

\subsection{Inequality constraints}\label{subsec:sim-inequality-constraints}

 In this section, we consider Problem~\eqref{eq:quadratic-online-optimization-general} with only inequality constraints, and we compare three different approaches to solve this optimization problem. 
 In the simulations we set $n = 10$ and we use the same matrix $\A$ used in Example~\ref{subsec:sim-equality-constraints} and as $\G^{\prime}$ we take the matrix $\G$ used in Example~\ref{subsec:sim-equality-constraints}. The signals $\bv_k$ and $\h^{\prime}_k$ are sinusoidal waves with a periodicity of $\omega=10^{-4}\pi$, with corresponding internal model $p\left(z\right) = z^2 - 2\cos\left(\omega\right)z + 1$.
Let $\x^*_{k,\mathrm{UNC}}$ be the optimal solution of the unconstrained version of~\eqref{eq:quadratic-online-optimization-general}, $\x^*_{k,\mathrm{UNC}} = \argmin_{\x \in \R^n} f_k(\x)$. The simulation is additionally set up so that periodically $\x_k^* = \x^*_{k,\mathrm{UNC}}$ and $\x_k^* \neq \x^*_{k,\mathrm{UNC}}$. 

\begin{figure}[!ht]
\centering
\includegraphics[scale = 0.55]{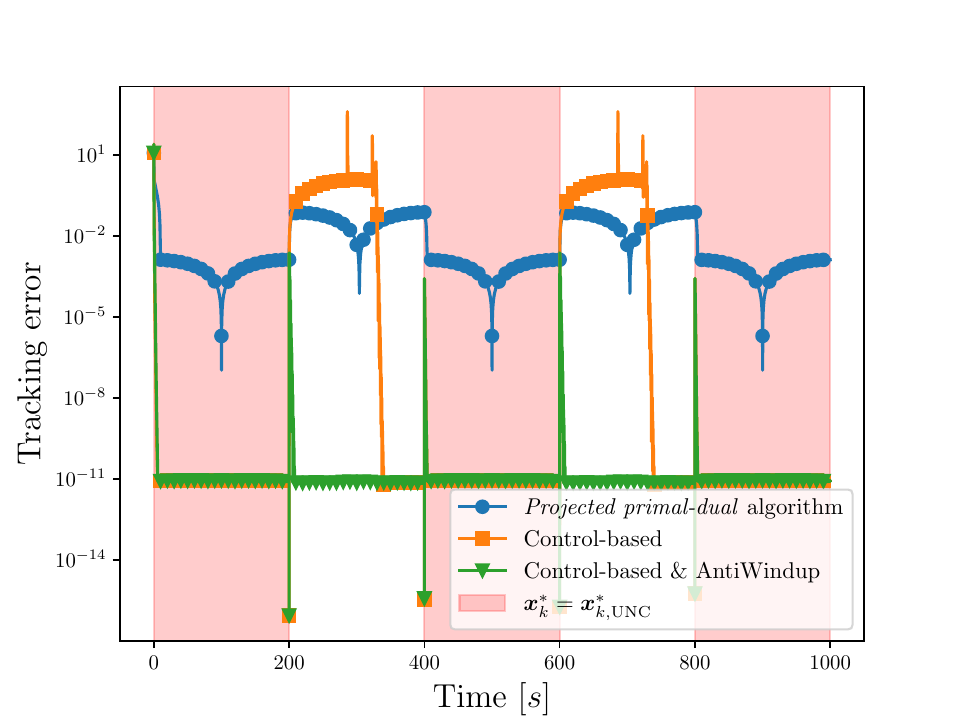}
\caption{Comparison between the \emph{projected primal-dual}, the control-based without anti wind-up, and the control-based with anti wind-up compensation in a semilogarithmic plot where $\bv_k$ and $\h^{\prime}_k$ are assumed to be two sine waves.}
\label{fig:comp-ineq}
\end{figure}

We compare three different algorithms:
\begin{itemize}
    \item Projected primal-dual algorithm~\eqref{eq:proPrimalDual}.
    \item Control-based algorithm modification of~\eqref{eq:control-based-algorithm}, which includes the saturation block but without anti-windup mechanism ($\rho = 0$).
    \item Control-based approach with anti wind-up compensation and $\rho=1$ presented in~\eqref{eq:control-based-algorithma-anti wind-up}.
\end{itemize}
Figure~\ref{fig:comp-ineq} depicts the tracking errors of these algorithms, with the shaded intervals highlighting the times when $\x^*_k=\x^*_{k,\mathrm{UNC}}$.

A first observation is that, similarly to the case with equality constraints, the (projected) primal-dual is outperformed by both control-based algorithms, in the sense that they can attain much smaller tracking errors. In particular, we can see that the control-based algorithms are converging to the optimal trajectory during the intervals of time when $\x^*_k=\x^*_{k,\mathrm{UNC}}$ and $\x^*_k \neq \x^*_{k,\mathrm{UNC}}$, while the unstructured primal-dual only reaches a bounded neighborhood of the trajectory. This again highlights the effectiveness of a control theoretical algorithm design.

Comparing the control-based algorithms, we notice that both undergo transients\footnote{We remark that the transients of ``Control-based'' in Figure~\ref{fig:comp-ineq} at times $400$, $800$ are present but hidden by the other trajectory.} when the optimization problem switches back and forth between the modes $\x^*_k=\x^*_{k,\mathrm{UNC}}$ and $\x^*_k \neq \x^*_{k,\mathrm{UNC}}$.
However, the algorithms differ significantly in the duration of these transients, as well as their magnitude. Indeed, the introduction of the anti-windup mechanism ensures much shorter transients, and smaller tracking errors.

In particular, the main issue when switching between $\x^*_k=\x^*_{k,\mathrm{UNC}}$ and $\x^*_k \neq \x^*_{k,\mathrm{UNC}}$ is that the dual variable $\w^{\prime}_k$ becomes saturated. In such cases, the algorithm continues to integrate the error $\f^{\prime}_k$ and increases the internal state variables associated with the integral actions. Ideally, when saturation is active, we are within the bounds of the inequality constraints, allowing the optimization algorithm to work as an unconstrained optimization. Consequently, the integration of the error $\f_k^{\prime}$ becomes not only redundant but also detrimental. The integral actions retain memory of the previously integrated error and negatively affect the algorithm's performance even outside the saturation region. The anti wind-up mechanism aims to reduce these integral actions when they are unnecessary, attempting to directly compensate for the input $\f_k^{\prime}$.

Finally, we observe that whenever $\x^*_k=\x^*_{k,\mathrm{UNC}}$, the control-based algorithms are such that $\w_k' = 0$, owing to the complementary slackness condition \cite{bertsekas_constrained_2014}. This implies that in this case the algorithms \eqref{eq:control-based-algorithm} and \eqref{eq:control-based-algorithma-anti wind-up} coincide, and in principle we can apply Proposition~\ref{pr:convergence-control-algorithm} to ensure exact convergence.
However, when the problems switches to $\x^*_k \neq \x^*_{k,\mathrm{UNC}}$ (and vice versa) a transient occurs. This implies that as long as the constraints are activated/deactivated throughout the online optimization, exact convergence to the optimal solution cannot be guaranteed.

\subsection{Time-varying quadratic term and constraint} 

As done in \cite{bastianello_internal_2022}, to assess the performance of our algorithm in more complex scenarios, we undertake an experimental validation. In these experiments we take $n = 10$ and $p = 2$. This begins by considering a time-varying quadratic term $\A_k = \A_1 + \sin\left(\omega k \right)\A_2$, where $\A_1$ is symmetric positive definite and with eigenvalues within the interval $\left[1,10\right]$, $\A_2$ is a symmetric sparse matrix of ones, where, the non zeros elements are in the ten percent of its entries. Moreover, we also consider time-varying constraints. Precisely, we set $\G_k = \G_1 + \sin\left(\omega k \right)\G_2$ where $\G_1$ is a random rectangular matrix with singular values in the range $\left[1,3\right]$ and, similar as before, $\G_2$ is a sparse matrix with the ten percent of its entries that are ones. The sinusoidal term $\sin\left(\omega k \right)$ has $\omega = 1/2$.
The remaining terms in Problem~\eqref{eq:quadratic-online-optimization} are held constant, namely $\bv_k = \bar{\bv}\in \R^n$ and $\h_k = \bar{\h}\in \R^p$. 

To implement our Algorithm we require an internal model that captures the dynamic nature of $f_k(\x)$. In view of the periodic evolution of $f_k(\x)$, reasonable choices of the internal model are
\begin{equation}\label{eq:models-timeVar}
    p\left(z\right) = \left(z-1\right)\prod_{l=1}^{L}\left[z^2 - 2\cos\left(l\omega\right)z+ 1 \right], \; L \in \mathbb{N}.
\end{equation}

Since $f_k(\x)$ is non-linear, we can only approximate the internal model. Utilizing model \eqref{eq:models-timeVar}, we aim to capture the model's minimum periodicity and fast variations.

\begin{figure}[!ht]
\centering
\includegraphics[scale = 0.55]{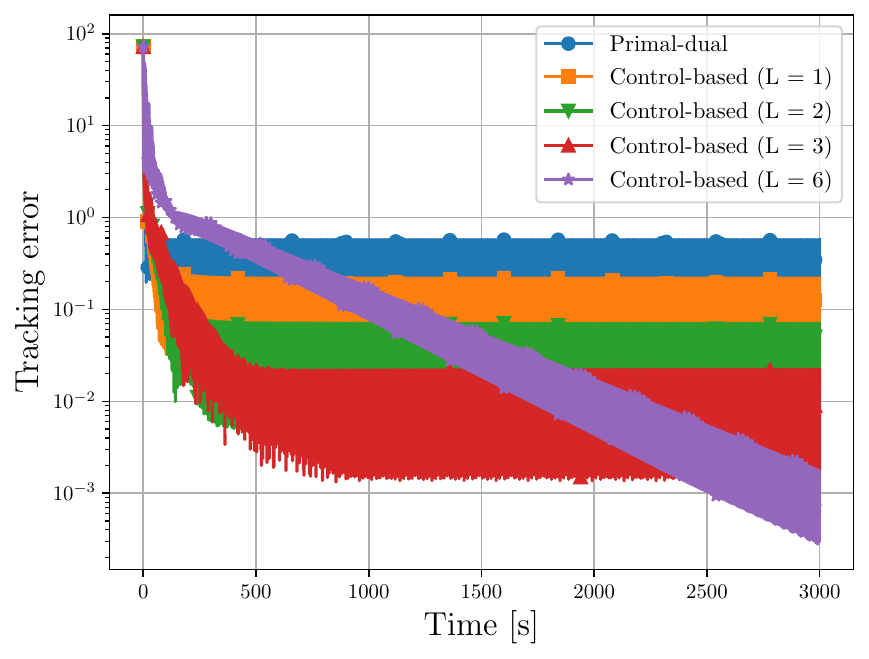}
\caption{Comparison between the \emph{primal-dual} algorithm, and the control-based with models of Equation~\eqref{eq:models-timeVar} ($L = 1,2,3,6$) in the time-varying case.}
\label{fig:comp-timeVar}
\end{figure}

The results of these simulations are depicted in Figure~\ref{fig:comp-timeVar}. As we can see, Algorithm~\eqref{eq:control-based-algorithm} exhibits better performance than the \emph{primal-dual} algorithm. Moreover, taking more complex internal models, we can improve the performance on the tracking error at the cost of a slower convergence rate. This slower convergence is characteristic of the internal model principle, as the controller requires a certain amount of time to reach its steady state response. This duration increases with the complexity of the internal model. Nevertheless, given our primary interest in the algorithm's asymptotic performances, this trade-off can be considered acceptable.

The proposed algorithm has good performance in the time-varying case (as well as the non-quadratic case of section~\ref{subsec:non-quadratic}). However, as highlighted in Remark~\ref{rem:genCase}, we have not been able to formally prove its convergence. This difficulty arises primarily because we approached the problem as a standard linear control problem, which does not hold in the time-varying scenario. Nonetheless, it may be possible to establish some bounds on the error by employing similar reasoning to the small gain approach discussed in ~\cite{bastianello_internal_2022}.

\subsection{Non-quadratic cost}\label{subsec:non-quadratic}
\begin{figure}[!ht]
\centering
\includegraphics[scale = 0.55]{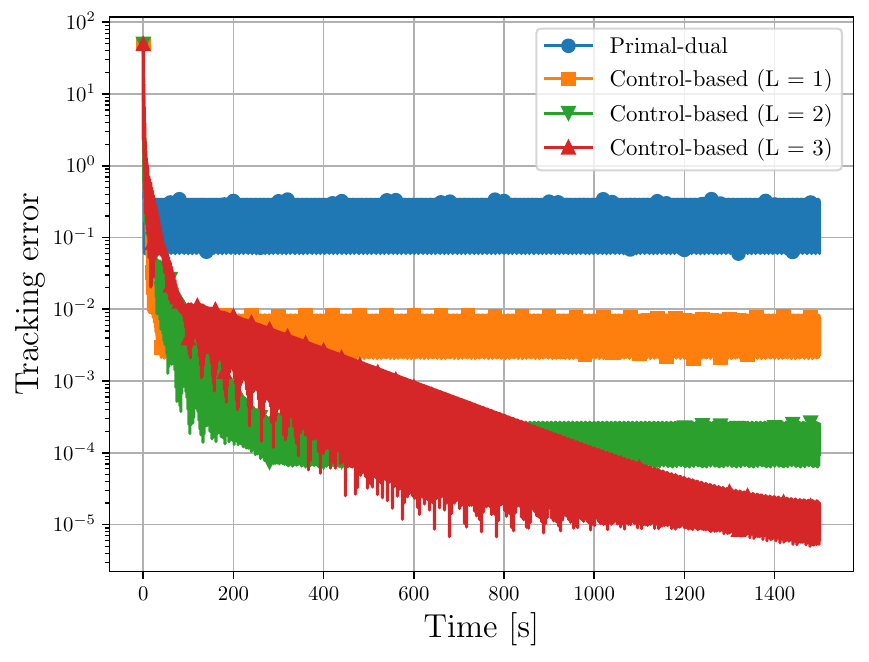}
\caption{Comparison between the \emph{ primal-dual} algorithm, and the control-based with models of Equation~\eqref{eq:models-timeVar} ($L = 1,2,3$) in the non-quadratic case.}
\label{fig:comp-nonQuad}
\end{figure}

In this section, we present an additional example of the application of our algorithm in a setting where the objective function is non-quadratic. Specifically, we utilize the following function taken from \cite{simonetto_class_2016}
\begin{equation*}
    f_k(\x) = \frac{1}{2}\x^{\top}\A\x + \bv^{\top}\x + \sin\left(\omega k\right)\log\left[1+\exp\left(\cv^{\top}\x\right)\right],
\end{equation*}
subject to the constraint:
$
    \G\x = \h_k.
$
In the above equations we take $\omega = 1/2$, $n = 10$, $\bv,\cv \in \R^{n}$ 
with $\cv$ such that $\norm{\cv}=1$. The matrix $\A \in \R^{n\times n}$ and $\G \in \R^{p\times n}$, with $p=1$, are generated as described in Section~\ref{subsec:sim-equality-constraints}, and the term $\h_k$ is the same as in case \ref{enum:case2-eq} of Section~\ref{subsec:sim-equality-constraints}. Finally, in Figure~\ref{fig:comp-nonQuad}, we report the tracking performance of the \emph{primal-dual} algorithm along with three different versions of the control-based algorithm that implement the internal models of Equation~\eqref{eq:models-timeVar}.
From Figure~\ref{fig:comp-nonQuad}, we can see that the control-based method outperforms the \emph{primal-dual} algorithm in terms of tracking performance, even if with a slower convergence rate.

\section{Conclusions}
In this paper we addressed the solution of online, constrained problems by leveraging control theory to design novel algorithms. When only equality constraints are present, we showed how to reformulate the problem as a robust, linear control problem, and we designed a suitable controller to achieve zero tracking error. When also inequality constraints are imposed, we showed how the required nonnegativity of the dual variables can lead to a wind-up phenomenon. As a consequence, we propose a modified version of the algorithm that incorporates an anti-windup scheme. Overall, we showed with both theoretical and numerical results how the proposed approaches outperform state-of-the-art alternatives.
The paper also points to interesting future developments. These include the application of our design principle to other classes of problems, and the exploration of alternative techniques to address inequality constraints.

\bibliographystyle{plain}
\bibliography{references}

\end{document}